\documentclass[psamsfonts]{amsart}
\usepackage{amssymb}
\usepackage{graphics}

\DeclareMathOperator{\Con}{Con}
\DeclareMathOperator{\ccon}{Con_c}
\DeclareMathOperator{\NId}{NId}
\DeclareMathOperator{\Id}{Id}
\DeclareMathOperator{\Idc}{Id_c}
\DeclareMathOperator{\FP}{FP}
\newcommand{\V}{\mathbf{V}}
\newcommand{\cC}{{\rm (C)}}

\newcommand{\urp}{Uniform Refinement Property}

\numberwithin{equation}{section}

\theoremstyle{plain}
\newtheorem{theorem}{Theorem}[section]
\newtheorem{proposition}[theorem]{Proposition}
\newtheorem{lemma}[theorem]{Lemma}
\newtheorem{corollary}[theorem]{Corollary}

\newtheorem*{SL}{Schmidt's Lemma}
\newtheorem*{thm}{Theorem}

\theoremstyle{definition}

\newtheorem{definition}[theorem]{Definition}

\newtheorem*{noteadd}{Note added}

\begin{document}

\title[Congruence lattices]%
{A uniform refinement property\\
for congruence lattices}

\author[F.~Wehrung]{Friedrich~Wehrung}
\address{D\'epartement de Math\'ematiques\\
         Universit\'e de Caen\\
         14032 Caen Cedex\\
         France}
 \email{gremlin@math.unicaen.fr}

 \date{\today}

\begin{abstract}
The Congruence Lattice Problem asks whether every
algebraic distributive lattice is isomorphic to
the congruence lattice of a lattice. It was
hoped that a positive solution would follow from
E. T. Schmidt's construction or from the approach
of P.~Pudl\'ak, M.~Tischendorf, and J.~T\r uma.
In a previous paper, we constructed a
distributive algebraic lattice $A$ with $\aleph_2$
compact elements that cannot be obtained by
Schmidt's construction.  In this paper, we show
that the same lattice $A$ cannot be obtained
using the Pudl\'ak, Tischendorf, T\r uma approach.

The basic idea is that every congruence lattice
arising from either method satisfies the
\emph{Uniform Refinement Property}, that is not
satisfied by our example. This yields, in
turn, corresponding negative results about
congruence lattices of sectionally complemented
lattices and two-sided ideals of von Neumann
regular rings.
 \end{abstract}

\maketitle

\section*{Introduction}
E. T. Schmidt introduces in \cite{Schm68} the
notion of a \emph{weakly distributive} (resp.,
\emph{distributive}) homomorphism of
semilattices, and proves the following important
result (see, for example, \cite[Satz 8.1]{Schm68}
and \cite[Theorem 3.6.9]{Schm81}):

\begin{SL} Let $A$ be an algebraic distributive
lattice. If the semilattice
of compact elements of $A$ is the image
under a distributive homomorphism of a
generalized Boolean semilattice, then there exists
a lattice $L$ such that \(\Con L\) is isomorphic
to $A$.
\end{SL}

This result yields important
partial positive answers to the Congruence
Lattice Problem (see, for example,
\cite{GrSc,Tisc94} for a survey). On the other
hand, we prove \cite[Theorem 2.15]{Wehr}, which
implies the following result:

\begin{thm} For any cardinal number
\(\kappa\geq\aleph_2\), there exists a
distributive semilattice \(S_\kappa\) of size
$\kappa$ that is not a weakly distributive image of
any distributive lattice.
\end{thm}

In this paper, we introduce a monoid-theoretical
property, the \emph{\urp} that \emph{is not}
satisfied by the semilattices \(S_\kappa\) of the
Theorem above, but it \emph{is} satisfied by
the congruence semilattice of any lattice
satisfying an additional condition, the
\emph{Congruence Splitting Property}
(Theorem~\ref{T:CSURP}). The latter is satisfied
by any lattice which is either sectionally
complemented, or relatively complemented, or a
direct limit of atomistic lattices
(Proposition~\ref{P:PptyCCS}). In particular, for
all \(\kappa\geq\aleph_2\), the
semilattice \(S_\kappa\) of the above Theorem
is not isomorphic to the congruence
semilattice of any sectionally complemented
lattice; this gives a partial negative answer to
\cite[Problem II.8]{Grat78}. It follows that no
semilattice \(S_\kappa\) can be realized as the
semilattice of finitely generated two-sided
ideals of a von Neumann regular ring.

\section*{Notation and terminology}

The \emph{refinement property} is the
monoid-theoretical axiom stating that for every
equation of the form
\(a_0+a_1=b_0+b_1\), there exist elements
\(c_{ij}\) ($i$, $j<2$) such that for all
\(i<2\), \(a_i=c_{i0}+c_{i1}\) and
\(b_i=c_{0i}+c_{1i}\). All semilattices will be
join-semilattices (not necessarily bounded). A
semilattice is, as usual, \emph{distributive},
if it satisfies the refinement property.

If $u$ and $v$ are elements of a lattice $L$,
then \(\Theta_L(u,v)\) (or \(\Theta(u,v)\) if
there is no ambiguity) denotes the least
congruence of $L$ identifying
$u$ and $v$. Furthermore, \(\Con L\) (resp.,
$\ccon L$) denotes the lattice (resp.,
semilattice) of all congruences (resp., compact
congruences) of $L$. If $L$ has a least element
(always denoted by $0$), an \emph{atom} of $L$ is
a minimal element of \(L\setminus\{0\}\).
A lattice $L$ with zero is \emph{atomistic}, if
every element of $L$ is a \emph{finite} join of
atoms.

Our lattice-theoretical results will be applied
to rings in Section~\ref{S:RegRings}. All our
rings are
\emph{associative} and \emph{unital} (but not
necessarily commutative). Recall that a ring $R$
is (von Neumann)
\emph{regular}, if it satisfies the axiom
\((\forall x)(\exists y)(xyx=x)\). If $R$ is
a regular ring, then every principal right ideal
of $R$ is of the form \(eR\), where $e$ is
idempotent, and the set
$\mathcal{L}(R)$ of all principal right ideals of
$R$, partially ordered by inclusion, is a
complemented modular lattice.

\section{Weak-distributive
homomorphisms}\label{S:wdhom}

We shall first modify slightly the original
definition, due to E. T. Schmidt
\cite{Schm81,Schm82} of a weakly distributive
homomorphism. A homomorphism of semilattices
\(f\colon S\to T\) is
\emph{weakly distributive} at an element $u$ of
$S$, if for all $y_0$, $y_1\in T$ such that
\(f(u)=y_0+y_1\), there are \(x_0\), \(x_1\in S\)
such that \(x_0+x_1=u\) and \(f(x_i)\leq y_i\),
for all \(i<2\). Say that $f$ is \emph{weakly
distributive}, if it is weakly distributive at
every element of $S$.

In particular, if $f$ is surjective, one
recovers the usual definition of a
weakly distributive homomorphism. Furthermore, note
that with this new definition, any composition of
two weakly distributive homomorphisms remains
weakly distributive.

\begin{lemma}\label{L:wdAdd} Let \(f\colon S\to
T\) be a homomorphism of semilattices, with $T$
distributive. Then the set of all elements of $S$
at which $f$ is weakly distributive is closed under
join.
\end{lemma}

\begin{proof} Let $X$ be the set of all elements
of $S$ at which $f$ is weakly distributive. It
suffices to prove that if $u'$ and $u''$ are any
two elements of $X$, then \(u=u'+u''\) belongs to
$X$. Thus let $y_i$ ($i<2$) be elements of $T$
such that \(f(u)=y_0+y_1\), that is, since
$f$ is a homomorphism of semilattices,
\(f(u')+f(u'')=y_0+y_1\) holds. Since
$T$ is distributive, it satisfies the refinement
property; thus there are decompositions
\(y_i=y'_i+y''_i\) ($i<2$) such that
\(f(u')=y'_0+y'_1\) and \(f(u'')=y''_0+y''_1\).
Since both $u'$ and $u''$ belong to $X$, there
are decompositions \(u'=x'_0+x'_1\) and
\(u''=x''_0+x''_1\) such that \(f(x'_i)\leq
y'_i\) and \(f(x''_i)\leq y''_i\) (\(i<2\)). Put
\(x_i=x'_i+x''_i\) (\(i<2\)). Then
\(f(x_i)\leq y_i\) and \(x_0+x_1=u\). Therefore,
\(u\in X\).
\end{proof}

The following result yields a large class of
weakly distributive homomorphisms:

\begin{proposition}\label{P:ConvHom} Let
\(e\colon K\to L\) be a lattice homomorphism with
convex range. Then the induced semilattice
homomorphism \(f\colon \ccon K\to\ccon L\) is
weakly dis\-tri\-bu\-ti\-ve.
\end{proposition}

\begin{proof} It is well-known that \(\ccon L\)
is a distributive semilattice (see, for example,
\cite[Theorem II.3.11]{Grat78}). Hence, by
Lemma~\ref{L:wdAdd}, it suffices to prove that
$f$ is weakly distributive at every
$\alpha$ of the form
\(\Theta_K(u,v)\) where \(u\leq v\) in $K$. Thus
let \(\beta_i\) ($i<2$) such that
\(f(\Theta_K(u,v))=\beta_0\vee\beta_1\),
that is,
\(\Theta_L(f(u),f(v))=\beta_0\vee\beta_1\). Thus,
by \cite[Lemma III.1.3]{Grat78}, there exist a
positive integer $n$ and elements \(w'_i\)
(\(i\leq 2n\)) of $L$ such that
\(f(u)=w'_0\leq w'_1\leq\cdots\leq w'_{2n}=f(v)\)
and for all \(i<n\),
\(w'_{2i}\equiv w'_{2i+1}\pmod{\beta_0}\) and
\(w'_{2i+1}\equiv w'_{2i+2}\pmod{\beta_1}\). But
the range of $e$ is convex in $L$, thus the
$w'_i$ belong to the range of $e$, thus there are
elements \(w_i\in K\) such that \(e(w_i)=w'_i\).
One can of course take \(w_0=u\) and
\(w_{2n}=v\), and, after replacing each \(w_i\) by
\(\bigvee_{j\leq i}(w_j\vee u)\wedge v\), one can
suppose that \(u=w_0\leq w_1\leq\cdots\leq
w_{2n}=v\). Then put
\(\alpha_0=\bigvee_{i<n}\Theta_K(w_{2i},w_{2i+1})\)
and
\(\alpha_1=\bigvee_{i<n}\Theta_K(w_{2i+1},w_{2i+2})\).
We have
\(\alpha_0\vee\alpha_1=\Theta_K(u,v)=\alpha\), and
\(f(\alpha_0)=\bigvee_{i<n}\Theta_L(f(w_{2i}),f(w_{2i+1}))
\subseteq\beta_0\); similarly,
\(f(\alpha_1)\subseteq\beta_1\).
\end{proof}

\section{The \urp}\label{S:URP}

In order to illustrate the terminology of the
section title, let us first consider any equation
system (in a given semilattice) of the form
\[
\Sigma\colon a_i+b_i=\mathrm{constant}
\qquad(\mbox{for all }i\in I).
\] When \(I=\{i,j\}\), a satisfactory notion of a
\emph{refinement} of $\Sigma$ consists of four
elements
\(c_{ij}^{uv}\) (\(u,v<2\)) satisfying the
equations
\begin{equation}
\begin{split}\label{Eq:Refab}
a_i=c_{ij}^{00}+c_{ij}^{01}\qquad&\mbox{and}\qquad
b_i=c_{ij}^{10}+c_{ij}^{11},\\
a_j=c_{ij}^{00}+c_{ij}^{10}\qquad&\mbox{and}\qquad
b_j=c_{ij}^{01}+c_{ij}^{11},
\end{split}
\end{equation}
 see Figure~1.
\begin{figure}[hbt]
\centerline{\includegraphics{Fig1.ill}}
\end{figure}

Note that (\ref{Eq:Refab}) implies immediately
the following consequence:
\begin{equation}\label{Eq:RefConseq} a_i\leq
a_j+c_{ij}^{01}.
\end{equation} When $I$ is an arbitrary finite
set, one can extend this in a natural way and
thus define a refinement of $\Sigma$ to be a
\(\mathcal{P}(I)\)-indexed family of elements of
$S$ satisfying suitable generalizations of
(\ref{Eq:Refab}). Nevertheless, this cannot be
extended immediately to the infinite case, so
that we shall focus instead on the consequence
(\ref{Eq:RefConseq}) of refinement, together with
an additional ``coherence condition"
\(c_{ik}^{01}\leq c_{ij}^{01}+c_{jk}^{01}\). More
precisely, we state the following definition:

\begin{definition}\label{D:URP} Let $S$ be a
semilattice, let $e$ be an element of $S$. Say
that the \emph{\urp} holds at $e$, if for all
families \((a_i)_{i\in I}\) and
\((b_i)_{i\in I}\) of elements of $S$ such that
\((\forall i\in I)(a_i+b_i=e)\), there are
families
\((a^*_i)_{i\in I}\),
\((b^*_i)_{i\in I}\) and 
\((c_{ij})_{(i,j)\in I\times I}\) of elements of
$S$ satisfying the following properties:
\begin{itemize}
\item[(i)] For all \(i\in I\), \(a^*_i\leq a_i\)
and
\(b^*_i\leq b_i\) and \(a^*_i+b^*_i=e\).

\item[(ii)] For all \(i,\ j\in I\),
\(c_{ij}\leq a^*_i,b^*_j\) and
\(a^*_i\leq a^*_j+c_{ij}\).

\item[(iii)] For all \(i,\ j,\ k\in I\),
\(c_{ik}\leq c_{ij}+c_{jk}\).
\end{itemize} Say that $S$ satisfies the \urp, if
the \urp\ holds at every element of $S$.
\end{definition}

It is to be noted that this is far from being the
only possible ``reasonable" definition for a
``\urp", see
\cite[Theorem 2.8, Claim 1]{Wehr}, which suggests
a quite different \urp\ (in the context of
partially ordered vector spaces). The following
result allows us to focus the investigation on a
generating set in the case where our semilattice
is distributive:

\begin{proposition}\label{P:URPAdd} Let $S$ be a
distributive semilattice. Then the set $X$ of all
elements of $S$ at which the \urp\ holds is
closed under join.
\end{proposition}

\begin{proof} Let $e_0$ and $e_1$ be two elements
of $X$, and put
\(e=e_0+e_1\). Let \((a_i)_{i\in I}\) and
\((b_i)_{i\in I}\) be two families of elements of
$S$ such that for all \(i\in I\),
\(a_i+b_i=e\). Since $S$ is distributive, there
are decompositions
\(a_i=a_{i0}+a_{i1}\), \(b_i=b_{i0}+b_{i1}\) (for
all
\(i\in I\)) such that
\(e_\nu=a_{i\nu}+b_{i\nu}\) (for all \(i\in I\)
and \(\nu<2\)). Since both \(e_0\) and \(e_1\)
belong to $X$, to the latter decompositions
correspond elements \(a^*_{i\nu}\),
\(b^*_{i\nu}\) and \(c_{ij\nu}\) (\(i,\ j\in I\)
and
\(\nu<2\)) witnessing the \urp\ at $e_0$ and
$e_1$.

Now put \(a^*_i=a^*_{i0}+a^*_{i1}\),
\(b^*_i=b^*_{i0}+b^*_{i1}\) and
\(c_{ij}=c_{ij0}+c_{ij1}\). It is obvious that
(i) to (iii) of Definition~\ref{D:URP} above are
satisfied with respect to \((a_i)_{i\in I}\) and
\((b_i)_{i\in I}\), thus proving that
\(X+X\subseteq X\).
\end{proof}

Note also the following result, whose easy proof
we shall omit:

\begin{proposition}\label{P:URPClwd} Let
\(f\colon S\to T\) be a weakly distributive
homomorphism of semilattices and let \(u\in S\).
If the \urp\ holds at $u$ in
$S$, then it also holds at \(f(u)\) in $T$.\qed
\end{proposition}

\section{Congruence splitting lattices; property
\cC}
\label{S:CSpptyC}

We will say throughout this paper that a lattice
$L$ is \emph{congruence splitting}, if for all
\(a\leq b\) in $L$ and all congruences
\(\alpha_0\) and
\(\alpha_1\) of $L$ such that
\(\Theta(a,b)\subseteq\alpha_0\vee\alpha_1\),
there are elements
$x_0$ and $x_1$ of \([a,\,b]\) such that
\(x_0\vee x_1=b\) and for all \(i<2\),
\(\Theta(a,x_i)\subseteq\alpha_i\). Of course, it
suffices to consider the case where both
\(\alpha_0\) and \(\alpha_1\) are compact
congruences and
\(\Theta(a,b)=\alpha_0\vee\alpha_1\).

For all elements $a$, $b$ and $c$ of a given
lattice, write
\(a\lessdot_cb\), if there exists $x$ such that
\(a\vee x=b\) and \(a\wedge x\leq c\).

\begin{definition} A lattice $L$ has
\emph{property
\cC}, if for all \(a\leq b\) and all $c$ in $L$,
there exist a positive integer $n$ and elements
\(x_i\) (\(i\leq n\)) of
$L$ such that
\(a=x_0\lessdot_cx_1\lessdot_c\cdots\lessdot_cx_n=b\).
\end{definition}

The letter `C' stands here for `complement'. In
the following proposition, we record a few
elementary properties of lattices either with
property \cC\ or congruence splitting.

\begin{proposition}\label{P:PptyCCS} The
following properties hold:
\begin{itemize}
\item[\rm (a)] Every relatively complemented
lattice, or every sectionally complemented
lattice, has property \cC.

\item[\rm (b)] Every atomistic lattice has
property \cC.

\item[\rm (c)] The class of all lattices
satisfying property
\cC\ is closed under direct limits.

\item[\rm (d)] Every lattice satisfying property
\cC\ is congruence splitting.

\item[\rm (e)] The class of all congruence
splitting lattices is closed under direct limits.
\end{itemize}
\end{proposition}

\begin{proof} (a) is obvious (take \(n=1\)).

(b) Let $L$ be an atomistic lattice. Note then
that for every \(a\in L\) and every atom $p$,
\(a\lessdot_0a\vee p\): indeed, if \(p\leq a\),
then in fact \(a=a\vee p\), and otherwise,
\(a\wedge p=0\). Now let \(a<b\) in $L$. There
exist a positive integer $n$ and atoms \(p_i\)
(\(i\leq n\)) such that
\(b=\bigvee_{i\leq n}p_i\). For all \(i\leq n\),
put
\(c_i=a\vee\bigvee_{j\leq i}p_j\). Then, by
previous remark,
\(a=c_0\lessdot_0c_1\lessdot_0\cdots\lessdot_0c_n=b\).

(c) is straightforward.

(d) Let $L$ be a lattice satisfying property \cC,
let $\alpha_j$ ($j<2$) be congruences of $L$, let
$a\leq b$ in $L$ such that
\(\Theta(a,b)\subseteq\alpha_0\vee\alpha_1\). To
reach the desired conclusion, we argue by
induction on the minimal length \(n=n(a,b)\) of a
chain
\(a=c_0\lessdot_a c_1\lessdot_a\cdots\lessdot_a
c_n=b\) such that for all \(i<n\), there exists
\(j<2\) such that
\(c_i\equiv c_{i+1}\pmod{\alpha_j}\) (such a chain
exists by \cite[Lemma III.1.3]{Grat78} and
property \cC). If
\(a<b\), then there exists $c$  such that
\(a\leq c\), \(n(a,c)<n(a,b)\), \(c\lessdot_ab\),
and, without loss of generality,
\(c\equiv b\pmod{\alpha_0}\). By induction
hypothesis, there are \(y_j\in[a,\,c]\) (\(j<2\))
such that $y_0 \vee y_1 = c$ and 
\(y_j\equiv a\pmod{\alpha_j}\). 
Furthermore, there
exists $z$ such that \(z\wedge c\leq a\) and
\(z\vee c=b\). Put \(x_0=y_0\vee z\) and
\(x_1=y_1\); then \(x_j\in[a,\,b]\),
\(x_j\equiv a\pmod{\alpha_j}\), and
\(x_0\vee x_1=b\).

(e) Let $L$ be a direct limit of a direct system
\((L_i,e_{ij})_{i\leq j\ \mathrm{in}\ I}\) of
lattices and lattice homomorphisms ($I$ is a
directed poset), with limiting maps \(e_i\colon
L_i\to L\) (for all \(i\in I\)). It is well-known
that \(\ccon L\) is then the direct limit of the
\(\ccon L_i\) with the corresponding transition
maps and limiting maps. Once this observation is
made, it is routine to verify that if all the
$L_i$ are congruence splitting, then so is~$L$.
\end{proof}

\begin{theorem}\label{T:CSURP} Let $L$ be a
congruence splitting lattice. Then \(\ccon L\)
satisfies the \urp.
\end{theorem}

\begin{proof} Put \(S=\ccon L\). By
Proposition~\ref{P:URPAdd}, it suffices to prove
that $S$ satisfies the \urp\ at every element of
$S$ of the form
\(\varepsilon=\Theta(u,v)\) where \(u\leq v\) in
$L$. Thus let
\((\alpha_i)_{i\in I}\) and
\((\beta_i)_{i\in I}\) be two families of elements
of $S$ such that for all \(i\in I\),
\(\alpha_i+\beta_i=\varepsilon\). Since $L$ is
congruence splitting, there are elements $s_i$
and $t_i$ of \([u,\,v]\) such that \(s_i\vee
t_i=v\) and
\(\Theta(u,s_i)\subseteq\alpha_i\) and
\(\Theta(u,t_i)\subseteq\beta_i\). Now, for all
\(i,j\in I\), put
\[
\alpha_i^*=\Theta(u,s_i),\qquad
\beta_i^*=\Theta(u,t_i)\qquad\mbox{and}\qquad
\gamma_{ij}=\Theta(s_j,s_i\vee s_j).
\] It is immediate that
\(\alpha_i^*\subseteq\alpha_i\),
\(\beta_i^*\subseteq\beta_i\) and
\(\alpha_i^*+\beta_i^*=\varepsilon\).
Furthermore, for all
\(i,j\in I\), we have
\(\gamma_{ij}\subseteq\Theta(u,s_i)=\alpha_i^*\)
and
\(\gamma_{ij}\subseteq\Theta(u,t_j)=\beta_j^*\).
Finally,
\(\gamma_{ij}\) is the least congruence $\theta$
of $L$ such that \(\theta(s_i)\leq\theta(s_j)\);
it follows immediately that
\(\gamma_{ik}\leq\gamma_{ij}+\gamma_{jk}\), thus
completing the proof.
\end{proof}

On the other hand, an inspection of the proof of
the Theorem stated in the Introduction
\cite[Theorem 2.15]{Wehr} shows in fact the
following property of the corresponding
semilattices \(S_\kappa\):

\begin{lemma}\label{L:SknonURP} For every
cardinal number \(\kappa\geq\aleph_2\), the
semilattice $S_\kappa$ does not satisfy the \urp\
at its largest element.\qed
\end{lemma}

By putting together Proposition~\ref{P:URPClwd}
and Lemma~\ref{L:SknonURP}, one deduces the
following:

\begin{corollary}\label{C:SknonMeas} For every
cardinal number \(\kappa\geq\aleph_2\), there are
no congruence splitting lattice $L$ and no
weakly distributive semilattice homomorphism
\(\mu\colon \ccon L\to S_\kappa\) with range
containing the largest element of $S_\kappa$.\qed
\end{corollary}

\begin{corollary}\label{C:Dichot} Let
\(\kappa\geq\aleph_2\) be a cardinal number.
Consider both following statements:
\begin{itemize}
\item[\rm (i)] There exist a lattice $L$ and a
weakly distributive homomorphism \(\mu\colon \ccon
L\to S_\kappa\) with range containing (as an
element) the largest element of
$S_\kappa$.

\item[\rm (ii)] For every bounded lattice $L$,
there exists a congruence splitting lattice $L'$
such that 
\(\ccon L\cong\ccon L'\).
\end{itemize} Then {\rm (i)} and {\rm (ii)}
cannot be simultaneously true.
\end{corollary}

\begin{proof} Suppose that both (i) and (ii) are
simultaneously true, and let $L$, $\mu$ be as in
(i). Since $L$ is the direct union of its closed
intervals, we obtain, by using
Proposition~\ref{P:ConvHom}, the existence of a
closed interval $K$ of $L$ such that the
restriction of $\mu$ to $K$ satisfies (i). Then,
applying (ii) to $K$ contradicts
Corollary~\ref{C:SknonMeas}.
\end{proof}

In particular, the Congruence Lattice Problem and
the problem whether every lattice has a
congruence-preserving embedding into a sectionally
complemented lattice cannot both have
positive answers.

\section{Applications to von Neumann regular
rings}
\label{S:RegRings}

For any ring $R$, we denote by \(\FP (R)\) the
class of all
\emph{finitely generated projective} right
$R$-modules, and by
\(\V(R)\) the (commutative) monoid of isomorphism
classes of elements of \(\FP (R)\), the addition
of \(\V(R)\) being defined by \([A]+[B]=[A\oplus
B]\) (\([A]\) denotes the isomorphism class of
$A$). Then \(\V(R)\) is \emph{conical},
that is, it satisfies the axiom
\((\forall x,y)(x+y=0\Rightarrow x=y=0)\).
Moreover, if $R$ is regular, then \(\V(R)\)
satisfies the refinement property. Furthermore,
we equip \(\V(R)\) with its \emph{algebraic}
preordering $\leq$, defined by \(x\leq y\) if and
only if there exists $z$ such that \(x+z=y\).
References about this can be found in
\cite{AGPO,Good79,Good95,VNeu60}.

\begin{lemma}\label{L:IsomTS} Let $R$ be a ring,
let \(J\in\Id R\), let $a$ and $b$ be idempotent
elements of $R$ such that \(aR\cong bR\) (as
right $R$-modules). Then \(a\in J\) if and only if
\(b\in J\).
\end{lemma}

\begin{proof} Since \(aR\cong bR\) and $a$ and
$b$ are idempotent, there are elements \(x\in
aRb\) and \(y\in bRa\) such that \(a=xy\) and
\(b=yx\). Suppose, for example, that
\(a\in J\). Then \(x\in aRb\subseteq J\) (because
$J$ is a right ideal of $R$), thus \(b=yx\in J\)
(because $J$ is a left ideal of $R$).
\end{proof}

If $L$ is any lattice with $0$, an ideal
\(\mathfrak{a}\) of $L$ is \emph{neutral}, if
\(x\in\mathfrak{a}\) and \(x\sim y\) implies that
\(y\in\mathfrak{a}\) (\(\sim\) is the relation of
perspectivity). We will denote by
\(\NId L\) the lattice of all neutral ideals of
$L$. Recall \cite{Birk93,Grat78} that if $L$ is a
sectionally complemented modular lattice, then
$\Con L$ and
$\NId L$ are (canonically) isomorphic.

\begin{lemma}\label{L:CharNeut} Let $R$ be a
regular ring. Then an ideal $\mathfrak{a}$ of
$\mathcal{L}(R)$ is neutral if and only if it is
closed under isomorphism (that is,
$J\in\mathfrak{a}$ and $I\cong J$ implies
$I\in\mathfrak{a}$).
\end{lemma}

\begin{proof} It is trivial that if
$\mathfrak{a}$ is closed under isomorphism, then
it is neutral. Conversely, suppose that
$\mathfrak{a}$ is neutral. Let
\(I,J\in\mathcal{L}(R)\) such that \(I\cong J\)
and
\(J\in\mathfrak{a}\). There exists
\(I'\in\mathcal{L}(R)\) such that \((I\cap
J)\oplus I'=I\). Since \(I\cap J\leq J\) and
\(J\in\mathfrak{a}\), we have
\(I\cap J\in\mathfrak{a}\). Furthermore, there
exists
\(J'\leq J\) such that \(I'\cong J'\). Since
\(I'\cap J=\{0\}\), we also have \(I'\cap
J'=\{0\}\), thus, by
\cite[Proposition 4.22]{Good79}, \(I'\sim J'\).
Since
\(J'\leq J\) and \(J\in\mathfrak{a}\), we have
\(J'\in\mathfrak{a}\), thus, since
\(\mathfrak{a}\) is neutral,
\(I'\in\mathfrak{a}\). Hence,
\(I=(I\cap J)\oplus I'\in\mathfrak{a}\).
\end{proof}

We denote by $\Id R$ the (algebraic) lattice of
two-sided ideals of any ring $R$, and by $\Idc R$
the semilattice of all compact (that is,
finitely generated) elements of
$\Id R$.

\begin{theorem}\label{T:NIdisId} Let $R$ be a
regular ring. Then one can define two mutually
inverse isomorphisms by the following rules:
\begin{gather*}
\varphi\colon \NId\mathcal{L}(R)\to\Id R,\qquad
\mathfrak{a}\mapsto\{x\in R\colon
xR\in\mathfrak{a}\},\\
\mbox{and}\\
\psi\colon \Id R\to\NId\mathcal{L}(R),\qquad
I\mapsto\{J\in\mathcal{L}(R)\colon J\subseteq I\}.
\end{gather*}
\end{theorem}

\begin{proof} First of all, we must verify that
for all \(\mathfrak{a}\in\NId\mathcal{L}(R)\),
\(\varphi(\mathfrak{a})\) as defined above is a
two-sided ideal of $R$. It is obvious that
\(\varphi(\mathfrak{a})\) is an additive subgroup
of $R$. Let
\(x\in\varphi(\mathfrak{a})\) and \(\lambda\in
R\). Then
\([x\lambda R]\leq[xR]\) and, since $R$ is
regular,
\([\lambda xR]\leq[xR]\) (the natural surjective
homomorphism
\(xR\to\lambda xR\) splits). Therefore, it is
sufficient to prove that if \(J\in\mathfrak{a}\)
and \([I]\leq[J]\), then
\(I\in\mathfrak{a}\). However, this results
immediately from Lemma~\ref{L:CharNeut}.
Conversely, the fact that $\psi$ takes its values
in \(\NId\mathcal{L}(R)\) results immediately
from Lemma~\ref{L:IsomTS}. The verification of
the fact that both $\varphi$ and $\psi$ are
order-preserving and mutually inverse is
straightforward.
\end{proof}

\begin{corollary}\label{C:ConcIdc} For any
regular ring $R$,
$\ccon \mathcal{L}(R)$ is isomorphic to $\Idc
R$.\qed
\end{corollary}

In \cite{Berg}, G. M. Bergman proves among other
things that any \emph{countable} bounded
distributive semilattice is isomorphic to $\Idc
R$ for some regular (and even ultramatricial)
ring $R$. By Corollaries \ref{C:SknonMeas} and
\ref{C:ConcIdc}, one cannot generalize this to
semilattices of size $\aleph_2$ (the size
$\aleph_1$ case is still open):

\begin{corollary}\label{C:NonRingMeas} For any
cardinal number \(\kappa\geq\aleph_2\), $S_\kappa$
is not isomorphic to the semilattice of finitely
generated two-sided ideals of any regular
ring.\qed
\end{corollary}

Note that this also provides a negative answer to
\cite[Problem 2.16]{Wehr}, \emph{via} the
following easy result:

\begin{proposition}\label{P:IdcQuotV} Let $R$ be
a regular ring. Then \(\Idc R\) is isomorphic to
the maximal semilattice quotient of \(\V(R)\).
\end{proposition}

\begin{proof} Let \(\pi:\V(R)\to\mathcal{P}(R)\)
be the mapping defined by the rule
\[
\pi(\alpha)=\{x\in R\colon  [xR]\leq
n\alpha\mbox{ for some positive integer }n\}.
\] Using the refinement property of \(\V(R)\), it
is easy to verify that $\pi$ is a monoid
homomorphism taking its values in \(\Id R\).
Moreover, for every
\(I\in\mathcal{L}(R)\), we have \(\pi([I])=RI\)
and those elements generate \(\Idc R\), thus $\pi$
maps \(\V(R)\) onto \(\Idc R\). Finally, by
\cite[Corollary 2.23]{Good79}, for all
\(I,\ J\in\mathcal{L}(R)\),
\(\pi([I])\leq\pi([J])\) if and only if there
exists a positive integer $n$ such that
\([I]\leq n[J]\); again by using refinement, it
follows that for all elements
\(\alpha,\ \beta\in\V(R)\),
\(\pi(\alpha)\leq\pi(\beta)\) if and only if
there exists a positive integer $n$ such that
\(\alpha\leq n\beta\). The conclusion follows.
\end{proof}

\begin{noteadd} Recently, the author, in a joint
paper with M. Plo\v s\v cica and J.
T\r uma, has proved that there exists a bounded
lattice $L$ such that \(\Con L\) is not
isomorphic to \(\Con L'\), for any congruence
splitting lattice $L'$.  Compare this with 
Corollary~\ref{C:Dichot} of this paper.
\end{noteadd}

\end{document}